\documentclass[12pt]{article}
\usepackage{amssymb,latexsym}
\input epsf

\title{\bf{Some combinatorial properties of flag simplicial pseudomanifolds
and spheres}}

\author{Christos~A.~Athanasiadis
\footnote{Supported by the 70/4/8755 ELKE Research Fund of the University
of Athens}\\
\small Department of Mathematics (Division of Algebra-Geometry)\\[-0.8ex]
\small University of Athens, Panepistimioupolis\\[-0.8ex]
\small 15784 Athens, Greece\\[-0.8ex]
\small \texttt{caath@math.uoa.gr}}

\date{\small May 27, 2009 }

\newcommand{\cC}{{\mathcal C}}

\newcommand{\gG}{{\mathcal G}}

\newcommand{\KK}{{\mathbb K}}

\newcommand{\sm}{{\smallsetminus}}
\newcommand{\qed}{$\hfill \Box$}

\begin{document}
\maketitle

\newtheorem{theorem}{Theorem}[section]
\newtheorem{proposition}[theorem]{Proposition}
\newtheorem{corollary}[theorem]{Corollary}
\newtheorem{definition}[theorem]{Definition}
\newtheorem{remark}[theorem]{Remark}
\newtheorem{lemma}[theorem]{Lemma}
\newtheorem{example}[theorem]{Example}
\newtheorem{examples}[theorem]{Examples}
\newtheorem{conjecture}[theorem]{Conjecture}
\newtheorem{question}[theorem]{Question}

\centerline{\small \emph{Dedicated to Anders Bj\"orner on the occasion
of his sixtieth birthday}}

\begin{abstract}
A simplicial complex $\Delta$ is called flag if all minimal nonfaces
of $\Delta$ have at most two elements. The following are proved: First,
if $\Delta$ is a flag simplicial pseudomanifold of dimension $d-1$, then
the graph of $\Delta$ (i) is $(2d-2)$-vertex-connected and (ii) has a 
subgraph which is a subdivision of the graph of the $d$-dimensional 
cross-polytope. Second, the $h$-vector of a flag simplicial homology 
sphere $\Delta$ of dimension $d-1$ is minimized when $\Delta$ is the 
boundary complex of the $d$-dimensional cross-polytope. 
\end{abstract}
%

%\maketitle

\section{Introduction}
\label{intro}

We will be interested in finite simplicial complexes. Such a complex
$\Delta$ is called \emph{flag} if every set of vertices which are 
pairwise joined by edges in $\Delta$ is a face of $\Delta$. For instance, 
every order complex (meaning the simplicial complex of all chains in a 
finite partially ordered set) is a flag complex. According to \cite[p. 
100]{Sta3}, flag complexes form a fascinating class of simplicial 
complexes which deserves further study. The class of flag complexes 
coincides with that of clique complexes of finite graphs.

Much of the combinatorial structure of flag complexes seems to be
significantly different from that of general simplicial complexes.
For instance, a simplicial sphere of dimension $d-1$ can
have as few as $d+1$ vertices and the minimum is attained by the
boundary complex of the $d$-dimensional simplex. In contrast, as
observed in \cite[Lemma 2.1.14]{Ga}, every flag simplicial (homology)
sphere (or, more generally, flag simplicial pseudomanifold) of
dimension $d-1$ has at least $2d$ vertices and the minimum is attained
by the boundary complex of the $d$-dimensional cross-polytope.
This paper proves some analogous (stronger) statements related to the
graph structure and face enumeration of flag complexes, which demonstrate
further that these complexes have a special position within the class
of all simplicial complexes.

We will denote by $\gG(\Delta)$ the one-dimensional skeleton
of a simplicial complex $\Delta$, called the \emph{graph} of $\Delta$.
Recall that, given a positive integer $m$, an abstract graph $\gG$ is
said to be $m$-connected if $\gG$ has at least $m+1$ nodes and any graph 
obtained from $\gG$ by deleting $m-1$ or fewer nodes and their incident 
edges is connected (necessarily with at least one edge).
It follows from Balinski's theorem \cite[Theorem 3.14]{Zi} that $\gG
(\Delta)$ is $d$-connected if $\Delta$ is the boundary complex of a
$d$-dimensional simplicial polytope and from \cite[Corollary 5]{Bar}
that the same statement holds for every simplicial pseudomanifold
$\Delta$ of dimension $d-1$. Our first result is the following.

\begin{theorem}
For every flag simplicial pseudomanifold $\Delta$ of dimension $d-1$,
the graph $\gG(\Delta)$ is $(2d-2)$-connected.
\label{thm1}
\end{theorem}

A \emph{subdivision} of an abstract graph $\gG$ is any graph which can
be obtained from $\gG$ by selecting some of the edges of $\gG$ and 
replacing each selected edge $e$ by a path with the same endpoints as $e$, 
so that the interiors of these paths are pairwise disjoint and do not 
intersect the set of nodes of $\gG$. Clearly, the graph of any simplicial 
pseudomanifold of dimension $d-1$ has a subgraph which is a
subdivision of the complete graph on $d+1$ nodes (the same property
was proved by Gr\"unbaum \cite{Gr} for graphs of $d$-dimensional convex
polytopes and by Barnette \cite{Bar} for a more general class of graphs
of cell decompositions of manifolds). Our second result is the following.

\begin{theorem}
For every flag simplicial pseudomanifold $\Delta$ of dimension $d-1$,
the graph $\gG(\Delta)$ has a subgraph which is a subdivision of the
graph of the $d$-dimensional cross-polytope.
\label{thm4}
\end{theorem}

It is not hard to show (see Proposition \ref{prop:LBT}) that for every
integer $k$, the number of faces of dimension $k$ of a
$(d-1)$-dimensional flag simplicial pseudomanifold $\Delta$ is minimized
when $\Delta$ is the boundary complex of the $d$-dimensional cross-polytope.
The same conclusion was proved by Meshulam \cite{Me} for the class of flag
simplicial complexes with nonzero reduced $(d-1)$-homology. Our third
result asserts that the analogous (stronger) statement for the $h$-vector
(see Section \ref{sec:pre} for definitions) of $\Delta$ is also valid when
one is restricted to a certain class of simplicial complexes which includes
all flag homology spheres, namely that of doubly Cohen-Macaulay flag
complexes. We refer to recent work of Nevo \cite{Ne} for conjectured lower
bounds on face numbers of flag (and more general) homology spheres with
fixed dimension and number of vertices.
\begin{theorem}
The h-vector $(h_0(\Delta), h_1(\Delta),\dots,h_d(\Delta))$ of any doubly
Cohen-Macaulay flag simplicial complex $\Delta$ of dimension $d-1$
satisfies the inequalities
\begin{equation}
h_i (\Delta) \, \ge \, {d \choose i}
\label{eq:2CM}
\end{equation}
for $0 \le i \le d$. In particular, these inequalities are valid
for all flag simplicial homology spheres of dimension $d-1$.
\label{thm2}
\end{theorem}

The previous theorem provides some (although only weak) evidence for 
the truth of a conjecture of Kalai (see \cite[p. 100]{Sta3}), stating 
that the $h$-vector of any flag Cohen-Macaulay simplicial complex is also 
equal to the $f$-vector of a (balanced) simplicial complex.

This paper is structured as follows. Section \ref{sec:pre} reviews
basic definitions and background on simplicial complexes, as well
as graph-theoretic terminology. Theorems \ref{thm1} and \ref{thm4} are
proved in Section \ref{sec:proof1}. Theorem \ref{thm2} is proved in
Section \ref{sec:proof2}. A higher dimensional analogue of Theorem 
\ref{thm1} is discussed in Section \ref{sec:high}.

\medskip
\noindent \textbf{Acknowledgements}. The author thanks Ronald 
Wotzlaw for useful discussions. The part of the proof of Theorem
\ref{thm1} which extends this result from the class of connected 
homology manifolds to that of pseudomanifolds is based on his ideas 
(see \cite{Wo}) on how the main result of \cite{Ath} can be extended 
to the setting of \cite{Bar}. The author also thanks Isabella Novik for 
the content of Remark \ref{rem:novik} and Ed Swartz and the anonymous
referee for useful comments.

\section{Preliminaries}
\label{sec:pre}

We will use the notation $[d] = \{1, 2,\dots,d\}$, when $d$ is a
positive integer, and write $|S|$ for the cardinality of a finite set
$S$.

\medskip
\textbf{Simplicial complexes}. Let $E$ be a finite set. An
(abstract) \emph{simplicial complex} on the ground set $E$ is a
collection $\Delta$ of subsets of $E$ such that $\sigma \subseteq
\tau \in \Delta$ implies $\sigma \in \Delta$. The elements of
$\Delta$ are called \emph{faces}. The dimension of a face $\sigma$
is defined as one less than the cardinality of $\sigma$. The
dimension of $\Delta$ is the maximum dimension of a face and is
denoted by $\dim (\Delta)$. Faces of $\Delta$ of dimension zero or
one are called \emph{vertices} or \emph{edges}, respectively. A
\emph{facet} of $\Delta$ is a face which is maximal with respect
to inclusion. The complex $\Delta$ is \emph{pure} if all its
facets have the same dimension. The $k$-skeleton $\Delta^{\le k}$
of $\Delta$ is the subcomplex formed by the faces of $\Delta$ of
dimension at most $k$. The simplicial join $\Delta_1 \ast
\Delta_2$ of two simplicial complexes $\Delta_1$ and $\Delta_2$ on
disjoint ground sets has as its faces the sets of the form
$\sigma_1 \cup \sigma_2$, where $\sigma_1 \in \Delta_1$ and
$\sigma_2 \in \Delta_2$.

The \emph{closed star} of $v \in E$ in $\Delta$ is the subcomplex
of $\Delta$ consisting of all subsets of those faces of $\Delta$
which contain $v$. The \emph{antistar} of $v$ in $\Delta$ is
defined as the restriction $\{\tau \in \Delta: v \notin \tau\}$ of
$\Delta$ on the set $E \sm \{v\}$ and is denoted by $\Delta \sm
v$. More generally, for $\sigma \subseteq E$ we denote by $\Delta
\sm \sigma$ the restriction $\{\tau \in \Delta: \tau \cap \sigma =
\varnothing\}$ of $\Delta$ on the set $E \sm \sigma$. The
\emph{link} of a face $\sigma$ in $\Delta$ is defined as $\Delta /
\sigma = \{\tau \sm \sigma: \tau \in \Delta, \, \sigma \subseteq
\tau\}$. For simplicity, we write $\Delta / v$ instead of $\Delta
/ \{v\}$ for $v \in E$.

A sequence $(\tau_0, \tau_1,\dots,\tau_n)$ of facets of $\Delta$ is
said to be a \emph{strong chain} if $\tau_{i-1} \cap \tau_i$ is a
codimension one face of both $\tau_{i-1}$ and $\tau_i$ for $1 \le i
\le n$. One can define an equivalence relation $\sim$ on the set of
facets of $\Delta$ by letting $\sigma \sim \tau$ if there exists a
strong chain $(\tau_0, \tau_1,\dots,\tau_n)$ of facets of $\Delta$
such that $\tau_0 = \sigma$ and $\tau_n = \tau$. The simplicial
complex formed by all subsets of the facets of $\Delta$ in an
equivalence class of $\sim$ is called a \emph{strong component} of
$\Delta$. A simplicial complex is \emph{strongly connected} if it has
a unique strong component. In particular, such a complex must be pure.
A strongly connected $(d-1)$-dimensional simplicial complex $\Delta$
is said to be a (simplicial) \emph{pseudomanifold} if each face of
$\Delta$ of dimension $d-2$ is contained in exactly two facets. It is
an easy observation that if $\Delta$ is a pseudomanifold, then $\Delta
\sm v$ is pure for every vertex $v$ of $\Delta$. The following stronger
statement is a special case of \cite[Lemma 2]{Bar}.
\begin{lemma} {\rm (cf. \cite[Lemma 2]{Bar})}
For every simplicial pseudomanifold $\Delta$ and vertex $v$, the complex
$\Delta \sm v$ is strongly connected.
\label{lem:anti}
\end{lemma}

The link $\Delta / \sigma$ of a face $\sigma$ in a pseudomanifold $\Delta$
may not be a pseudomanifold, since it may fail to be strongly connected.
However, any strong component of such a link is also a pseudomanifold.

The \emph{boundary complex} of a simplicial convex polytope $P$ is the 
abstract simplicial complex on the vertex set of $P$ for which a subset 
$\sigma$ of the vertex set of $P$ is a face if and only if the convex 
hull of $\sigma$ is a face of $P$, other than $P$ itself. For instance, 
the boundary complex of a $d$-dimensional simplex consists of all subsets 
of a $(d+1)$-element set of cardinality at most $d$. When we talk about
topological properties of an abstract simplicial complex $\Delta$,
we implicitly refer to those of its geometric realization
$\|\Delta\|$ \cite[Section 9]{Bj1}. For instance, $\Delta$ is said
to be a $(d-1)$-\emph{sphere} if $\|\Delta\|$ is homeomorphic to a
sphere of dimension $d-1$. We call $\Delta$ a \emph{homology
sphere} (over some fixed field $\KK$) if for all $\sigma \in
\Delta$ (including $\sigma = \varnothing$) we have
\[ \widetilde{H}_i \, (\Delta / \sigma, \KK) \ = \
\cases{ 0, & if \ $i < \dim \, (\Delta / \sigma)$, \cr
\KK, & if \ $i = \dim \, (\Delta / \sigma)$, } \]
where $\widetilde{H}_* (\Gamma, \KK)$ denotes reduced simplicial homology
of $\Gamma$ with coefficients in the field $\KK$. We call $\Delta$ a
\emph{homology manifold} (over $\KK$) if $\Delta / \sigma$ is a homology
sphere for every nonempty face $\sigma$ of $\Delta$. We note that every
homology manifold which is connected and has at least two vertices must
be a pseudomanifold.

Suppose that the ground set $E$ of $\Delta$ has $n$ elements, say $E =
\{v_1, v_2,\dots,v_n\}$. The \emph{face ring} (or \emph{Stanley-Reisner
ring}) associated to $\Delta$ is defined as the quotient $\KK [\Delta] =
\KK [x_1, x_2,\dots,x_n] / I_\Delta$ of the polynomial ring over $\KK$ in
the variables $x_1, x_2,\dots,x_n$ by the ideal $I_\Delta$ generated by
the square-free monomials $x_{i_1} x_{i_2} \cdots x_{i_r}$ for which
$\{i_1, i_2,\dots,i_r\} \notin \Delta$. The complex $\Delta$ is said to
be \emph{Cohen-Macaulay} or \emph{Gorenstein} (over $\KK$) if $\KK [\Delta]$
is a \emph{Cohen-Macaulay} or \emph{Gorenstein} ring, respectively. We
refer to \cite{Sta3} for a thorough discussion of these concepts. By
Reisner's theorem \cite[Corollary II.4.2]{Sta3}, $\Delta$ is Cohen-Macaulay
if and only if $\widetilde{H}_i \, (\Delta / \sigma, \KK) = 0$ for all
$\sigma \in \Delta$ and $i < \dim \, (\Delta / \sigma)$. Given a positive
integer $m$, the complex $\Delta$ is said to be
\emph{$m$-Cohen-Macaulay} (or \emph{doubly Cohen-Macaulay}, for $m=2$)
if $\Delta \sm \sigma$ is Cohen-Macaulay of the same dimension as $\Delta$
for all subsets $\sigma$ of $E$ of cardinality less than $m$ (including
$\sigma = \varnothing$). We have the hierarchy of properties

\begin{center}
\begin{tabular}{c}
sphere \, $\Rightarrow$ \, homology sphere \, $\Rightarrow$ \,
$\cases{
\begin{tabular}{l}
doubly Cohen-Macaulay \, $\Rightarrow$ \, Cohen-Macaulay  \\
homology manifold \, $\Rightarrow^*$ \, pseudomanifold 
\end{tabular} } $ \\ \\
\, $\Rightarrow$ \, pure
\end{tabular}
\end{center}

\noindent (where the implication $\Rightarrow^*$ assumes
connectivity and positive dimension). Moreover, boundary complexes
of simplicial polytopes are spheres. The classes of boundary
complexes of simplicial polytopes, homology spheres, homology
manifolds, $m$-Cohen-Macaulay complexes (for any fixed $m \ge 1$)
and pure complexes are all closed under taking links of faces. The
class of homology spheres coincides with that of nonacyclic
Gorenstein complexes.

Let $\Delta$ be any simplicial complex of dimension $d-1$. The number of
$k$-dimensional faces of $\Delta$ will be denoted by $f_k (\Delta)$, so
that $f_{-1} (\Delta) = 1$ unless $\Delta$ is the void complex
$\varnothing$. The sequence
$f (\Delta) = (f_{-1} (\Delta), f_0 (\Delta),\dots,f_{d-1} (\Delta))$
is called the $f$-\emph{vector} of $\Delta$. The $h$-\emph{vector} of
$\Delta$ is the sequence $h (\Delta) = (h_0 (\Delta), h_1
(\Delta),\dots,h_d (\Delta))$ defined by the equality
\begin{equation}
\sum_{i=0}^d \ h_i (\Delta) \, x^i \ = \
\sum_{i=0}^d \ f_{i-1} (\Delta) \, x^i (1-x)^{d-i}.
\label{eq:h-vec}
\end{equation}
The \emph{reduced Euler characteristic} of $\Delta$ is defined as
\begin{equation}
\widetilde{\chi} (\Delta) \ = \ (-1)^{d-1} \, h_d (\Delta) \ = \
\sum_{i=0}^d \ (-1)^{i-1} f_{i-1} (\Delta).
\label{eq:euler}
\end{equation}
The polynomial which appears in either hand-side of (\ref{eq:h-vec}) is
called the $h$-\emph{polynomial} of $\Delta$ and is denoted by $h_\Delta
(x)$. It is a fundamental property of Cohen-Macaulay complexes that
$h_i (\Delta) \ge 0$ holds for every index $i$.

A simplicial complex $\Delta$ is called \emph{flag} if all its minimal
nonfaces have at most two elements. The simplicial join of two flag
complexes and the link of any face of a flag complex are also flag
complexes. For instance, the boundary complex of the $d$-dimensional
cross-polytope \cite[Example 0.4]{Zi} is isomorphic to the simplicial
join of $d$ copies of the 0-sphere (the zero-dimensional complex
with just two vertices) and hence it is a flag complex. The next
proposition asserts that among all $(d-1)$-dimensional flag simplicial
pseudomanifolds, the boundary complex of the $d$-dimensional
cross-polytope has the minimum number of faces in each dimension.
\begin{proposition}
Any $(d-1)$-dimensional flag simplicial pseudomanifold has no fewer 
than $2^i {d \choose i}$ faces of dimension $i-1$ for all $0 \le i \le d$.
\label{prop:LBT}
\end{proposition}
\emph{Proof.} Let $\Delta$ be a flag simplicial pseudomanifold of dimension
$d-1$. We proceed by induction on $d$. The result is easily verified for
$i=0$ or $d=1$, so we assume that $i \ge 1$ and $d \ge 2$. Since $\Delta$
is a flag simplicial complex which is
not a simplex, there exist two vertices, say $u$ and $v$, of $\Delta$
such that $\{u, v\}$ is not an edge of $\Delta$. Among the
$(i-1)$-dimensional faces of $\Delta$, there exist $f_{i-2} (\Delta /
u)$ faces which contain $u$, $f_{i-2} (\Delta / v)$ faces which contain
$v$ and $f_{i-1} (\Delta / u)$ faces which belong to $\Delta / u$.
Since these three sets of $(i-1)$-dimensional faces of $\Delta$ are
pairwise disjoint, we conclude that
\[ f_{i-1} (\Delta) \ \ge \ f_{i-1} (\Delta / u) \, + \, f_{i-2}
(\Delta / u) \, + \, f_{i-2} (\Delta / v). \]
Since every strong component of $\Delta / u$ or $\Delta / v$ is a flag
pseudomanifold of dimension $d-2$, it follows from the previous
inequality and the induction hypothesis that
\begin{eqnarray*}
f_{i-1} (\Delta) & \ge & 2^i {d-1 \choose i} \, + \, 2 \cdot 2^{i-1}
{d-1 \choose i-1} \\
& & \\
&=& 2^i {d \choose i}.
\end{eqnarray*}
This completes the induction and the proof.
\qed

\medskip
%\noindent
\textbf{Graphs}.
A graph is a simplicial complex of dimension zero or one. We will 
refer to the vertices of a graph as \emph{nodes}, to avoid
possible confusion with vertices of other simplicial complexes
considered simultaneously. Two nodes $u$ and $v$ of a graph $\gG$ are
said to be \emph{adjacent} (or joined by an edge) in $\gG$ if
$\{u, v\}$ is an edge of $\gG$. A \emph{walk} of length
$n$ in $\gG$ is an alternating sequence $w = (v_0, e_1,
v_1,\dots,e_n, v_n)$ of nodes and edges, such that $e_i = \{v_{i-1},
v_i\}$ for $1 \le i \le n$. We say that $w$ \emph{connects} nodes $v_0$
and $v_n$, which are the \emph{endpoints} of $w$. The walk $w$ is said
to be a \emph{path} if $v_0, v_1,\dots,v_n$ are pairwise distinct; in
this case $v_1,\dots,v_{n-1}$ are the \emph{interior nodes} of $w$. We say
that $\gG$ is \emph{connected} if any two nodes can be connected by a
walk in $\gG$. Given a positive integer $m$, the graph $\gG$ is said to
be \emph{$m$-connected} if it has at least $m+1$ nodes and $\gG \sm
\sigma$ is connected for all subsets $\sigma$ of the set of nodes of
$\gG$ with cardinality less than $m$. Equivalently, $\gG$ is
$m$-connected if it is one-dimensional and $m$-Cohen-Macaulay over some
field (equivalently, over all fields) as a simplicial complex. A 
\emph{subgraph} of $\gG$ is any graph which can be obtained by deleting 
some of the edges of $\gG \sm \sigma$, for some subset $\sigma$ of the 
set of nodes of $\gG$. A \emph{subdivision} of $\gG$ is any graph which 
can be obtained from $\gG$ by selecting some of the edges of $\gG$ and 
replacing each selected edge $e$ by a path with the same endpoints as $e$, 
so that the interiors of these paths are pairwise disjoint and do not 
intersect the set of nodes of $\gG$. 

\medskip
%\noindent
\textbf{The graph $\gG(\Delta)$}. The 1-skeleton of a simplicial
complex $\Delta$ is called the \emph{graph} of $\Delta$ and is
denoted by $\gG(\Delta)$. We are primarily interested in
$\gG(\Delta)$ when $\Delta$ is a pseudomanifold. A walk $(v_0,
e_1, v_1,\dots,e_n, v_n)$ in $\gG(\Delta)$ will be called a
\emph{$\Delta$-strong walk} if for every index $1 \le i \le n$
there exists a facet $\tau_i$ of $\Delta$ containg $e_i$, so that
$\tau_i$ and $\tau_{i+1}$ lie in the same strong component of the
closed star of $v_i$ in $\Delta$ for all $1 \le i \le n-1$. This
condition allows for arguments which use induction on the
dimension of a pseudomanifold, by considering the links of its
vertices. It implies, for instance, that for every index $1 \le i
\le n-1$, some strong component of $\Delta / v_i$ contains both
$v_{i-1}$ and $v_{i+1}$.

Part (i) of the following proposition is implicit in the proof of
\cite[Theorem 4]{Bar}. We include the proof for the convenience of
the reader.

\begin{proposition}
Let $\Delta$ be a simplicial pseudomanifold of dimension $d-1$.
\begin{itemize}
\itemsep=0pt
\item[{\rm (i)}] {\rm (cf. \cite[Theorem 4]{Bar})}
If $\sigma$ is any subset of the vertex set of $\Delta$ of cardinality
less than $d$, then any two vertices of $\Delta$ not in $\sigma$ can be
connected by a $\Delta$-strong walk in $\gG(\Delta) \sm \sigma$. In
particular, the graph $\gG (\Delta)$ is $d$-connected.

\item[{\rm (ii)}]
If $\sigma$ is a face of $\Delta$, then any two vertices of $\Delta$
not in $\sigma$ can be connected by a $\Delta$-strong walk in $\gG (\Delta) 
\sm \sigma$. In particular, the graph $\gG(\Delta) \sm \sigma$ is 
connected.
\end{itemize}
\label{prop:ba}
\end{proposition}
\emph{Proof.} Let $a$ and $b$ be any two vertices of $\Delta$ not in
$\sigma$. We will show that $a$ and $b$ can be connected by a
$\Delta$-strong walk in $\gG (\Delta) \sm \sigma$. Pick any element
$v$ of $\sigma$. By Lemma \ref{lem:anti}, there exists a strong
chain $\cC = (\tau_0, \tau_1,\dots,\tau_n)$ of facets of $\Delta \sm
v$ such that $a \in \tau_0$ and $b \in \tau_n$. We claim that the
intersection $\tau_{i-1} \cap \tau_i$ is not a subset of $\sigma$ for
any index $1 \le i \le n$. Indeed, this is clear in part (i) since
$\tau_{i-1} \cap \tau_i$ has $d-1$ elements and does not contain $v$.
If $\sigma$ is a face of $\Delta$, as in part (ii), of cardinality at
least $d$, so that $\sigma$ is a facet, then
the claim holds because $\sigma$ does not appear in the chain $\cC$ and
$\tau_{i-1}$ and $\tau_i$ are the only facets of $\Delta$ which contain
$\tau_{i-1} \cap \tau_i$. For $1 \le i \le n$ we pick any vertex $v_i \in
\tau_{i-1} \cap \tau_i$ not contained in $\sigma$ and observe that $\{a,
v_1, \dots,v_n, b\}$ is the set of nodes of a $\Delta$-strong walk in
$\gG (\Delta) \sm \sigma$ connecting $a$ and $b$. This completes the
proof.
\qed

\section{Connectivity of $\gG(\Delta)$}
\label{sec:proof1}

In this section we prove Theorems \ref{thm1} and \ref{thm4}.

\medskip
\noindent \emph{Proof of Theorem \ref{thm1}.} We will first prove the
theorem for connected homology manifolds and then indicate how this proof
can be modified in the case of pseudomanifolds.

Let $\Delta$ be a connected
flag simplicial homology manifold of dimension $d-1$ and let $\tau$ be
any subset of the vertex set of $\Delta$ of cardinality less than $2d-2$.
Since $\Delta$ has at least $2d$ vertices, we only need to show that
any two vertices, say $a$ and $b$, of $\Delta$ not in $\tau$ can be
connected by a walk in $\gG(\Delta) \sm \tau$.
We proceed by induction on $d$. Clearly, we may assume that $d \ge 3$.
Let $\sigma$ denote the set of elements of $\tau$ which are adjacent to at
least $2d-4$ elements of $\tau$ in $\gG(\Delta)$. Observe that, in view
of our assumption on the cardinality of $\tau$, each element of $\sigma$
is adjacent in $\gG(\Delta)$ to all other elements of $\tau$. In 
particular, the elements of $\sigma$ are pairwise adjacent in 
$\gG(\Delta)$. As a result, since $\Delta$ is a flag complex, $\sigma$ 
a face of $\Delta$. Therefore, by Proposition \ref{prop:ba} (ii), any 
two vertices of $\Delta$ not in $\sigma$ can be connected by a walk in 
$\gG(\Delta) \sm \sigma$. Let $w = (v_0, e_1, v_1,\dots,e_n, v_n)$  be 
such a walk connecting $a$ and $b$. 

We may
assume that no two consecutive nodes of $w$ are in $\tau$. Indeed,
suppose that $v_{i-1}$ and $v_i$ are both elements of $\tau$ for some
index $i$. The link of $e_i$ in $\Delta$ is a flag homology sphere of
dimension $d-3$ and hence it has at least $2d-4$ vertices. Since $\tau$
has at most $2d-5$ elements other than $v_{i-1}$ and $v_i$, there exists
a vertex $u$ of $\Delta / e_i$ not in $\tau$. Then $\{u, v_{i-1},
v_i\}$ is a two-dimensional face of $\Delta$ and inserting $u$ between
$v_{i-1}$ and $v_i$ in $w$ results in a walk in $\gG(\Delta) \sm \sigma$
having a smaller number of pairs of consecutive nodes in $\tau$.
Repeating this process for every pair of consecutive nodes of $w$ in
$\tau$ results in a walk in $\gG(\Delta) \sm \sigma$ connecting $a$ and
$b$ with the desired property.

Consider any node $v_i$ of $w$ which is an element of $\tau$, so that
$1 \le i \le n-1$ and neither $v_{i-1}$ nor $v_{i+1}$ is an element of
$\tau$. To complete the proof, it suffices to show that $v_{i-1}$ and
$v_{i+1}$ can be connected by a walk in $\gG(\Delta) \sm \tau$ for any
such index $i$. Since $\Delta / v_i$ is a flag simplicial homology sphere
of dimension $d-2$, by our induction hypothesis the graph $\gG(\Delta /
v_i)$ is $(2d-4)$-connected. By construction $v_i$ is not in $\sigma$
and hence at most $2d-5$ vertices of $\Delta / v_i$ are in $\tau$. As
a result, $v_{i-1}$ and $v_{i+1}$ can be connected by a walk in
$\gG(\Delta / v_i) \sm \tau$ and hence by a walk in $\gG(\Delta) \sm
\tau$.

Finally, suppose that $\Delta$ is a flag simplicial pseudomanifold
of dimension $d-1$. By Proposition \ref{prop:ba} (ii), we may choose
the walk $w$ in the previous argument to be $\Delta$-strong. Let
$\tau_1,\dots,\tau_n$ be facets of $\Delta$ as in the definition of a
$\Delta$-strong walk. When inserting a vertex $u$ not in $\tau$ between
two consecutive nodes $v_{i-1}$ and $v_i$ of $w$ which are in $\tau$,
we can guarantee that the new walk will also be $\Delta$-strong.
Indeed, let $\Gamma$ be the strong component of $\Delta / e_i$ which
contains $\tau_i \sm e_i$. Since $\Gamma$ is a flag pseudomanifold of
dimension $d-3$, it has at least $2d-4$ vertices and we may choose $u$
to be in $\Gamma$. By construction, $\{u, v_{i-1}, v_i\}$ is contained
in a facet of $\Delta$ which can be connected to $\tau_i$ by a strong
chain of facets, each of which contains both $v_{i-1}$ and $v_i$.
From this fact it follows that our new walk is also $\Delta$-strong. In
the final part of the argument we only need to replace the link of $v_i$
in $\Delta$ with its strong component which contains $\tau_i \sm 
\{v_i\}$ and $\tau_{i+1} \sm \{v_i\}$.
\qed

\medskip
The next statement follows from the proof of Theorem \ref{thm1}.

\begin{corollary}
Let $\Delta$ be a flag simplicial pseudomanifold of dimension $d-1$. If
$\tau$ is a subset of the vertex set of $\Delta$ of cardinality less than
$2d-2$, then any two vertices of $\Delta$ not in $\tau$ can be connected
by a $\Delta$-strong walk in $\gG (\Delta) \sm \tau$. \qed
\label{cor:thm1}
\end{corollary}
\begin{remark} {\rm
Responding to a question posed by the author in a previous version of this 
paper, Novik \cite{No} has shown that the $k$-skeleton of every 
$(d-1)$-dimensional flag simplicial homology sphere is 
$2(d-k)$-Cohen-Macaulay (this statement generalizes Theorem \ref{thm1} in 
the case of flag homology spheres). The proof uses the Stanley-Reisner ring 
of $\Delta$, \cite[Lemma 5.1]{NS} and the Taylor resolution for quadratic 
monomial ideals.
\qed}
\label{rem:novik}
\end{remark}

%\medskip
\noindent \emph{Proof of Theorem \ref{thm4}.}
Let $\sigma = \{v_1, v_2,\dots,v_d\}$ be a facet of $\Delta$. Since
$\Delta$ is a pseudomanifold, for each $1 \le i \le d$ there exists a
unique vertex $u_i$ of $\Delta$ other than $v_i$ such that $(\sigma \sm
\{v_i\}) \cup \{u_i\}$ is also a facet of $\Delta$. Since $\Delta$ is flag,
the $u_i$ are pairwise distinct. We set $\tau = \{u_1, u_2,\dots,u_d\}$
and recall that the graph of the $d$-dimensional cross-polytope can be
obtained from the complete graph on $2d$ nodes by removing $d$ edges which
are mutually disjoint. Since we have $\{v_i, v_j\} \in \Delta$ and $\{v_i, 
u_j\} \in \Delta$ for distinct indices $i$ and $j$, it suffices to show 
that any two elements of $\tau$ can be connected by a path in $\gG(\Delta)$ 
so that the sets of interior nodes of all ${d \choose 2}$ resulting paths 
are mutually disjoint and each such set intersects neither $\sigma$ nor 
$\tau$.

Consider the face $\sigma_{ij} = \sigma \sm \{v_i, v_j\}$ of $\Delta$,
where $1 \le i < j \le d$. The link $\Delta / \sigma_{ij}$ is a
one-dimensional simplicial complex each vertex of which belongs to
exactly two edges. As a result, $\Delta / \sigma_{ij}$ is a disjoint
union of one-dimensional spheres. Moreover, it contains the edges $\{v_i,
v_j\}$, $\{v_i, u_j\}$ and $\{v_j, u_i\}$. From these facts it follows that
there exists a path $p_{ij}$ in $\gG(\Delta)$ connecting $u_i$ and $u_j$,
each interior node of which is a vertex of $\Delta / \sigma_{ij}$ other
than $v_i$ and $v_j$. 

We claim that (i) the set of nodes of $p_{ij}$ does
not intersect $\sigma$, (ii) the set of interior nodes of $p_{ij}$ does
not intersect $\tau$ and (iii) no vertex of $\Delta$ is an interior node
of two or more of the paths $p_{ij}$. Indeed, (i) is clear by
construction. Consider $u_r \in \tau \sm \{u_i, u_j\}$. Since $\Delta$ is
$(d-1)$-dimensional and flag and $(\sigma \sm \{v_r\}) \cup \{u_r\} \in
\Delta$, we have $\{u_r, v_r\} \notin \Delta$. Since $v_r \in \sigma_{ij}$,
we conclude that $u_r$ is not a vertex of $\Delta / \sigma_{ij}$ and hence
that $u_r$ cannot be one of the nodes of $p_{ij}$. This proves (ii). Finally,
suppose that $u$ is a vertex of $\Delta$ which is an interior node of two
paths $p_{ij}$ and $p_{k \ell}$. Then $u$ belongs to both links $\Delta /
\sigma_{ij}$ and $\Delta / \sigma_{k \ell}$ and hence we have $\{u, v\} \in 
\Delta$ for every $v \in \sigma_{ij} \cup \sigma_{k \ell}$. Since $\Delta$ 
is flag, it follows that $\sigma_{ij} \cup \sigma_{k \ell} \cup \{u\} \in 
\Delta$. Since $\sigma_{ij} \cup \sigma_{k \ell}$ is equal to $\sigma$, if 
$\{i, j\}$ and $\{k, \ell\}$ are disjoint, and to $\sigma \sm \{v_r\}$ for 
some $r \in \{i, j, k, \ell\}$ otherwise, it follows that $\sigma \cup \{u\} 
\in \Delta$ or $(\sigma \sm \{v_r\}) \cup \{u\} \in \Delta$. By our choice 
of $u_r$ and since $\sigma$ is a facet of $\Delta$, we conclude that either
$u \in \sigma$ or $u = u_r$, contradicting (i) and (ii). This contradiction 
proves (iii). It follows from facts (i)-(iii) that the paths $p_{ij}$ have 
the desired properties.
\qed

\section{Lower bound for the $h$-vector}
\label{sec:proof2}

In this section we prove Theorem \ref{thm2}. We will make use of
the following elementary lemma.
\begin{lemma}
Let $\Delta$ be a pure simplicial complex and $v$ be a vertex of 
$\Delta$. We have
\[ h_\Delta (x) \ = \
\cases{ h_{\Delta \sm v} (x) + x \, h_{\Delta / v} (x), & if \
$\dim (\Delta \sm v) = \dim (\Delta)$, \cr
h_{\Delta \sm v} (x), & otherwise. } \]
\label{lem:del}
\end{lemma}
\emph{Proof.} Let $d-1 = \dim (\Delta)$, so that $\dim (\Delta / v) =
d-2$ and either $\dim (\Delta \sm v) = d-1$ or else $\Delta \sm v =
\Delta / v$. By considering those $(i-1)$-dimensional faces of $\Delta$
which contain $v$ and those which do not, we see that
\[ f_{i-1} (\Delta) \ = \ f_{i-1} (\Delta \sm v) + f_{i-2} (\Delta / v) \]
for $0 \le i \le d$ (where $f_{i-1} (\Gamma) = 0$ for negative integers
$i$ by convention). In either case, multiplying this equation with $x^i
(1-x)^{d-i}$, summing and using (\ref{eq:h-vec}) we arrive at the proposed
equality expressing the $h$-polynomial of $\Delta$ in terms of those of
$\Delta \sm v$ and $\Delta / v$.
\qed

\medskip
\noindent \emph{Proof of Theorem \ref{thm2}.} Let $\Delta$ be a doubly
Cohen-Macaulay flag simplicial complex of dimension $d-1$. We need to
show that
\begin{equation}
h_\Delta (x) \, \ge \, (1+x)^d,
\label{eq:tLBT}
\end{equation}
where such an inequality will be meant to hold coefficientwise. We proceed
by induction on $d$. The statement
holds for $d=1$ since then $\Delta$ consists of $q \ge 2$ vertices,
having no other nonempty faces, and $h_\Delta (x) = 1 + (q-1) x$.
Suppose that $d \ge 2$. Since $\Delta$ is a flag simplicial complex
which is not a simplex, there exist two vertices, say $u$ and $v$, of
$\Delta$ such that $\{u, v\}$ is not an edge of $\Delta$. Since the link
$\Delta / v$ is doubly Cohen-Macaulay and flag of dimension $d-2$, our
induction hypothesis implies that
\begin{equation}
h_{\Delta / v} (x) \, \ge \, (1+x)^{d-1}.
\label{eq:link}
\end{equation}
For the same reason we have
\[ h_{\Delta / u} (x) \, \ge \, (1+x)^{d-1}. \]
Let $\Gamma$ denote the closed star of $u$ in $\Delta$. Then
$\Gamma$ is a subcomplex of $\Delta \sm v$ and both $\Gamma$ and
$\Delta \sm v$ are Cohen-Macaulay of dimension $d-1$. The
monotonicity property of $h$-vectors \cite[Theorem 2.1]{Sta2}
implies that $h_{\Delta \sm v} (x) \ge h_\Gamma (x)$ and Lemma
\ref{lem:del} implies that $h_\Gamma (x) = h_{\Gamma \sm u} (x)  =
h_{\Delta / u} (x)$, so that
\begin{equation}
h_{\Delta \sm v} (x) \, \ge \, h_\Gamma (x) \ = \ h_{\Delta / u} (x) \,
\ge \, (1+x)^{d-1}.
\label{eq:del}
\end{equation}
The desired inequality (\ref{eq:tLBT}) follows by combining (\ref{eq:link}) 
and (\ref{eq:del}) with Lemma \ref{lem:del}.
\qed

\medskip
We conjecture that if equality holds in (\ref{eq:2CM}) for some $1 \le i 
\le d-1$, then $\Delta$ is isomorphic to the boundary complex of the 
$d$-dimensional cross-polytope. This statement does not follow immediately 
from the previous proof.

\section{A higher dimensional analogue}
\label{sec:high}

Balinski's theorem on the one-dimensional skeleton $\gG (P)$ of a
convex polytope $P$ was generalized in \cite{Ath} to the graphs $\gG_k
(P)$ defined as follows. The nodes of $\gG_k (P)$ are the $k$-dimensional
faces of $P$ and two such faces are adjacent if there exists a
$(k+1)$-dimensional face of $P$ which contains them both. Theorem
\ref{thm1} can also be generalized in this direction. Given any
$(d-1)$-dimensional simplicial complex $\Delta$ and integer $0 \le k
\le d-2$, we denote by $\gG_k (\Delta)$ the graph with nodes the
$k$-dimensional faces of $\Delta$, in which two such faces are
adjacent if there exists a $(k+1)$-dimensional face of $\Delta$ which
contains them both. The graphs $\gG_k (P)$ and $\gG_k (\Delta)$ reduce
to $\gG (P)$ and $\gG (\Delta)$, respectively, for $k=0$.
\begin{theorem}
If $n_k (d) = 2(k+1)(d-k-1)$, then the graph $\gG_k (\Delta)$ is $n_k
(d)$-connected for every connected flag simplicial homology manifold
$\Delta$ of dimension $d-1$ and all integers $0 \le k \le d-2$.
\label{thm3}
\end{theorem}

The proof of Theorem \ref{thm3} is similar to that of the main result 
of \cite{Ath} and is omitted. The value of $n_k (d)$ in Theorem 
\ref{thm3} cannot be improved, as the example of the $d$-dimensional 
cross-polytope shows. It is likely that Theorem \ref{thm3} can be 
extended to the class of all flag simplicial pseudomanifolds.

\end{document}